\documentclass[a4,12pt]{amsart}

\usepackage{amsfonts,amssymb,amscd,amsmath,enumerate,verbatim,calc}
\usepackage{amsmath}
\usepackage{amssymb}
\usepackage{amscd}
\usepackage{ascmac}

\newtheorem{Theorem}{Theorem}[section]
\newtheorem{Lemma}[Theorem]{Lemma}

\newtheorem{Remark}[Theorem]{Remark}

\newtheorem{Example}[Theorem]{Example}
\newtheorem{Claim}[Theorem]{Claim}

\newcommand{\superarrow}[2]{\mathrel{\mathop{#1}\limits_{\scriptstyle {#2}}}}

\def\RMN#1{\uppercase\expandafter{\romannumeral#1}}

\makeatletter
  
  \@addtoreset{equation}{section}
\makeatother

\renewcommand{\qed}{{\unskip\nobreak\hfil\penalty50\quad\null\nobreak\hfil{\bf
q.e.d.}\parfillskip0pt\finalhyphendemerits0\par\medskip}}

\newcommand{\spec}{{\rm Spec}}

\title{The divisor class groups and the graded canonical modules of multi-section rings}
\author{Kazuhiko Kurano}

\begin{document}

\maketitle

\renewcommand{\thefootnote}{\fnsymbol{footnote}}
\footnote[0]{

This work was supported by KAKENHI 21540050.}

\begin{abstract}
We shall describe the divisor class group and the graded canonical module of the multi-section ring $T(X; D_1, \ldots, D_s)$ 
defined in (\ref{T}) below for a normal projective variety $X$
and Weil divisors $D_1$, \ldots, $D_s$ on $X$
under a mild condition.
In the proof, we use the theory of Krull domain and
the equivariant twisted inverse functor due to Hashimoto~\cite{H}. 
\end{abstract}

\section{Introduction}
We shall describe the divisor class groups and the graded canonical modules of multi-section rings
associated with a normal projective variety.

Suppose that ${\mathbb Z}$, ${\mathbb N}_0$ and ${\mathbb N}$ 
are the set of integers, non-negative integers and positive integers,
respectively.

Let $X$ be a normal projective variety over a field $k$ with the function field $k(X)$.
We always assume $\dim X > 0$.
We denote by $C^1(X)$ the set of closed  subvarieties of $X$ of codimension $1$.
For $V \in C^1(X)$ and $a \in k(X)^\times$, we define as
\begin{eqnarray*}
{\rm ord}_V(a) & = & \ell_{{\mathcal O}_{X, V}}({\mathcal O}_{X, V}/\alpha{\mathcal O}_{X, V})
-\ell_{{\mathcal O}_{X, V}}({\mathcal O}_{X, V}/\beta{\mathcal O}_{X, V}) \\
{\rm div}_X(a) & = & \sum_{V \in C^1(X)}{\rm ord}_V(a) \cdot V  \in  {\rm Div}(X) = \bigoplus_{V \in C^1(X)}
{\Bbb Z}\cdot V ,
\end{eqnarray*}
where $\alpha$ and $\beta$ are elements in ${\mathcal O}_{X, V}$ such that 
$a = \alpha/\beta$, and
$\ell_{{\mathcal O}_{X, V}}( \ )$ denotes the length as an 
${\mathcal O}_{X, V}$-module.

We call an element in ${\rm Div}(X)$ a {\em Weil divisor} on $X$.
For a Weil divisor $D = \sum n_V V$, we say that 
$D$ is {\em effective}, and write $D \ge 0$, if $n_V \ge 0$ for any $V \in  C^1(X)$.
For a Weil divisor $D$ on $X$, we put
\[
H^0(X, {\mathcal O}_X(D)) = 
\{ a \in k(X)^\times \mid {\rm div}_X(a) + D \ge 0 \} \cup \{ 0 \} .
\]
Here we remark that $H^0(X, {\mathcal O}_X(D))$ is a 
$k$-vector subspace of $k(X)$.

Let $D_1$, \ldots, $D_s$ be Weil divisors on $X$.
We define the multi-section rings $T(X; D_1, \ldots, D_s)$ and $R(X; D_1, \ldots, D_s)$ associated with $D_1$, \ldots, $D_s$
as follows:
\begin{eqnarray}
& & T(X; D_1, \ldots, D_s) \label{T} \\
\nonumber & = & 
\bigoplus_{(n_1, \ldots, n_s) \in {{\mathbb N}_0}^s} H^0(X, {\mathcal O}_X(\sum_i n_iD_i))t_1^{n_1}\cdots t_s^{n_s}
\subset k(X)[t_1, \ldots, t_s]  \\
& & \nonumber \\
& & R(X; D_1, \ldots, D_s) \nonumber  \\ & = & 
\bigoplus_{(n_1, \ldots, n_s) \in {\mathbb Z}^s} H^0(X, {\mathcal O}_X(\sum_i n_iD_i))t_1^{n_1}\cdots t_s^{n_s}
\subset k(X)[t_1^{\pm 1}, \ldots, t_s^{\pm 1}] \nonumber 
\end{eqnarray}
We want to describe the divisor class groups and the graded canonical modules 
of the above rings.

For a Weil divisor $F$ on $X$, we set
\[
M_F = 
\bigoplus_{(n_1, \ldots, n_s) \in {\mathbb Z}^s} H^0(X, {\mathcal O}_X(\sum_i n_iD_i +F))
t_1^{n_1}\cdots t_s^{n_s}
\subset k(X)[t_1^{\pm 1}, \ldots, t_s^{\pm 1}] ,
\]
that is, $M_F$ is a ${\mathbb Z}^s$-graded 
reflexive $R(X; D_1, \ldots, D_s)$-module with
\[
[M_F]_{(n_1, \ldots, n_s)} = H^0(X, {\mathcal O}_X(\sum_i n_iD_i +F))
t_1^{n_1}\cdots t_s^{n_s} .
\]
We denote by $\overline{M_F}$ the isomorphism class of 
the reflexive module $M_F$ 
in ${\rm Cl}(R(X; D_1, \ldots, D_s))$.

For a normal variety $X$, we denote by 
${\rm Cl}(X)$ the class group of $X$, and
for a Weil divisor $F$ on $X$,
we denote by $\overline{F}$ the residue class represented by 
the Weil divisor $F$ 
in ${\rm Cl}(X)$.

In the case where ${\rm Cl}(X)$ is freely generated by $\overline{D_1}$, \ldots, $\overline{D_s}$,
the ring $R(X; D_1, \ldots, D_s)$ is usually called the {\em Cox ring} of $X$
and denoted by ${\rm Cox}(X)$.

\begin{Remark}\label{T(X; D)}
\begin{rm}
Assume that $D$ is an ample divisor on $X$.
In this case, $T(X; D)$ coincides with $R(X; D)$, and it is a Noetherian normal domain
by a famous result of Zariski (see Lemma~2.8 in \cite{HuKeel}).
It is well-known that ${\rm Cl}(T(X; D))$ is isomorphic to ${\rm Cl}(X)/{\Bbb Z}\overline{D}$.
Mori \cite{M} constructed a lot of examples of non-Cohen Macaulay factorial domains
using this isomorphism.

It is well-known that the canonical module of $T(X; D)$ is isomorphic to $M_{K_X}$,
and the canonical sheaf $\omega_X$ coincides with $\widetilde{M_{K_X}}$.
Watanabe proved a more general result in Theorem~(2.8) in \cite{W}.
\end{rm}
\end{Remark}

We want to establish  the same type of  the above results 
for multi-section rings.

For $R(X; D_1, \ldots, D_s)$, we had already proven the following:

\begin{Theorem}[Elizondo-Kurano-Watanabe~\cite{EKW}, Hashimoto-Kurano~\cite{HK}]
\label{EKWHK}
Let $X$ be a normal projective variety over a field such that $\dim X > 0$.
Assume that $D_1$, \ldots, $D_s$ are Weil divisors on $X$ 
such that ${\mathbb Z}D_1 + \cdots + {\mathbb Z}D_s$ contains an ample Cartier divisor.
Then, we have the following:
\begin{enumerate}
\item
$R(X; D_1, \ldots, D_s)$ is a Krull domain.
\item
The set $\{ P_V \mid V \in C^1(X) \}$ coincides with the set of homogeneous prime ideals
of $R(X; D_1, \ldots, D_s)$ of height $1$, where $P_V = M_{-V}$.
\item
We have an exact sequence
\[
0 \longrightarrow \sum_i {\mathbb Z}\overline{D_i}
\longrightarrow {\rm Cl}(X) \stackrel{p}{\longrightarrow} {\rm Cl}(R(X; D_1, \ldots, D_s))
\longrightarrow 0
\]
such that $p(\overline{F}) = \overline{M_F}$.
\item
Assume that $R(X; D_1, \ldots, D_s)$ is Noetherian.
Then $\omega_{R(X; D_1, \ldots, D_s)}$ is isomorphic to $M_{K_X}$ as a ${\mathbb Z}^s$-graded module.
Therefore, $\omega_{R(X; D_1, \ldots, D_s)}$ is $R(X; D_1, \ldots, D_s)$-free
if and only if $\overline{K_X} \in \sum_i {\mathbb Z}\overline{D_i}$ in ${\rm Cl}(X)$.
\end{enumerate}
\end{Theorem}

Suppose that ${\rm Cl}(X)$ is finitely generated free
${\Bbb Z}$-module generated by $\overline{D_1}$, \ldots, $\overline{D_s}$.
By the above theorem,
the Cox ring ${\rm Cox}(X)$ is factorial and 
\[
\omega_{{\rm Cox}(X)}  = M_{K_X} = {\rm Cox}(X) (\overline{K_X}),
\]
where we regard ${\rm Cox}(X)$ as a ${\rm Cl}(X)$-graded ring.

The main result of this paper is the following:

\begin{Theorem} \label{Main}
Let $X$ be a normal projective variety over a field $k$
such that $d = \dim X > 0$.
Assume that $D_1$, \ldots, $D_s$ are Weil divisors on $X$ 
such that ${\mathbb N}D_1 + \cdots + {\mathbb N}D_s$ contains an ample Cartier divisor.
Put
\[
U = \{ j \mid {\rm tr. deg}_kT(X; D_1, \ldots, D_{j-1}, D_{j+1}, \ldots, D_s) = d + s - 1 \} .
\]
Then, we have the following:
\begin{enumerate}
\item
$T(X; D_1, \ldots, D_s)$ is a Krull domain.
\item
The set 
\[
\{ Q_V \mid V \in C^1(X) \} \cup \{ Q_j \mid j \in U \}
\]
coincides with the set of homogeneous prime ideals
of $T(X; D_1, \ldots, D_s)$ of height $1$, where 
\[
Q_V = P_V \cap T(X; D_1, \ldots, D_s)
\]
and 
\[
Q_j  = 
\bigoplus_{\superarrow{n_1, \ldots, n_s \in {\Bbb N}_0}{n_j > 0}} 
T(X; D_1, \ldots, D_s)_{(n_1, \ldots, n_s)} .
\]
\item
We have an exact sequence
\[
0 \longrightarrow \sum_{j \not\in U} {\mathbb Z}\overline{D_j}
\longrightarrow {\rm Cl}(X) \stackrel{q}{\longrightarrow} {\rm Cl}(T(X; D_1, \ldots, D_s))
\longrightarrow 0
\]
such that $q(\overline{F}) = \overline{M_F \cap k(X)[t_1, \ldots, t_s, \{ t_j^{-1} \mid j \not\in U \} ]}$.
\item
Assume that $T(X; D_1, \ldots, D_s)$ is Noetherian.
Then $\omega_{T(X; D_1, \ldots, D_s)}$ is isomorphic to 
\[
M_{K_X} \cap t_1 \cdots t_sk(X)[t_1, \ldots, t_s, \{ t_j^{-1} \mid j \not\in U \} ]
\]
as a ${\mathbb Z}^s$-graded module.
Further, we have 
\[
q(\overline{K_X + \sum_iD_i}) = \overline{\omega_{T(X; D_1, \ldots, D_s)}} .
\]
Therefore, $\omega_{T(X; D_1, \ldots, D_s)}$ is $T(X; D_1, \ldots, D_s)$-free
if and only if 
\[
\overline{K_X + \sum_iD_i}
\in \sum_{j \not\in U} {\mathbb Z}\overline{D_j}
\]
 in ${\rm Cl}(X)$.
\end{enumerate}
\end{Theorem}

Here, ${\rm tr. deg}_kT$ denotes the transcendence degree
of the fractional field of $T$ over a field $k$.

\begin{Remark}
\begin{rm}
With notation as in the previous theorem,
${\rm ht}(Q_j) = 1$ if and only if $j \in U$.
This will be proven in Lemma~\ref{Le}.
Since ${\mathbb N}D_1 + \cdots + {\mathbb N}D_s$ contains an ample Cartier divisor, $Q_j \neq (0)$ for any $j$.
Therefore, 
${\rm ht}(Q_j) \ge 2$ if and only if $j \not\in U$.
\end{rm}
\end{Remark}

\section{Examples}

\begin{Example}\label{ex1}
\begin{rm}
Let $X$ be a normal projective variety with $\dim X > 0$.
Assume that all of $D_i$'s are ample Cartier divisors on $X$.
Then,  $T(X; D_1, \ldots, D_s)$ is Noetherian by a famous result of Zariski
(see Lemma~2.8 in \cite{HuKeel}).

Assume that $s = 1$.
By definition, $U = \emptyset$ since $\dim X > 0$.
By Theorem~\ref{Main} (3), ${\rm Cl}(T(X; D_1))$ is isomorphic to ${\rm Cl}(X)/{\Bbb Z}\overline{D_1}$.
By Theorem~\ref{Main} (4), $\omega_{T(X; D_1)}$ is $T(X; D_1)$-free module 
if and only if 
\[
\overline{K_X} \in {\Bbb Z}\overline{D_1}
\]
in ${\rm Cl}(X)$ (see Remark~\ref{T(X; D)}).

Next, assume that $s \ge 2$.
In this case, $U = \{ 1, 2, \ldots, s \}$.
By Theorem~\ref{Main} (3), ${\rm Cl}(X)$ is isomorphic to ${\rm Cl}(T(X; D_1, \ldots, D_s))$.
By Theorem~\ref{Main} (4), $\omega_{T(X; D_1, \ldots, D_s)}$ is $T(X; D_1, \ldots, D_s)$-free module 
if and only if 
\[
\overline{K_X} =  \overline{- D_1 - \cdots - D_s}
\]
in ${\rm Cl}(X)$.
When this is the case, $-K_X$ is ample, that is, $X$ is a Fano variety.
\end{rm}
\end{Example}

\begin{Example}\label{ex2}
\begin{rm}
Set $X = {\Bbb P}^m \times {\Bbb P}^n$.
Let $p_1$ (resp.\  $p_2$) be the first (resp.\ second) projection.

Let $H_1$ be a hyperplane of ${\Bbb P}^m$, and $H_2$ a hyperplane of ${\Bbb P}^n$.
Put $A_i = p_i^{-1}(H_i)$ for $i = 1, 2$.
In this case, ${\rm Cl}(X) = {\Bbb Z} \overline{A_1} + {\Bbb Z} \overline{A_2} \simeq {\Bbb Z}^2$,
and $K_X = -(m+1) A_1 - (n+1) A_2$.

We have
\[
{\rm Cox}(X) = R(X; A_1, A_2) = k[x_0, x_1, \ldots, x_m, y_0, y_1, \ldots, y_n] .
\]
${\rm Cox}(X)$ is a ${\Bbb Z}^2$-graded ring such that 
$x_i$'s (resp.\ $y_j$'s) are of degree $(1,0)$ (resp.\ $(0,1)$).

Let $a$, $b$, $c$, $d$ be positive integers such that $ad - bc \neq 0$.
Put $D_1 = aA_1 + bA_2$ and $D_2 = cA_1 + dA_2$.
Then, both $D_1$ and $D_2$ are ample divisors.
Consider the multi-section rings:
\begin{eqnarray*}
R(X; D_1, D_2) & = & \oplus_{p, q \in {\Bbb Z}} {\rm Cox}(X)_{p(a,b) + q(c,d)} \\
T(X; D_1, D_2) & = & \oplus_{p, q \ge 0} {\rm Cox}(X)_{p(a,b) + q(c,d)} 
\end{eqnarray*}
Here, both $R(X; D_1, D_2)$ and $T(X; D_1, D_2)$ are Cohen-Macaulay rings.

By Theorem~\ref{EKWHK} (4), we know
\begin{eqnarray*}
\mbox{$R(X; D_1, D_2)$ is a Gorenstein ring} & \Longleftrightarrow &
\mbox{$\overline{K_X} \in {\Bbb Z}\overline{D_1} + {\Bbb Z}\overline{D_2}$ in ${\rm Cl}(X)$} \\
& \Longleftrightarrow &
(m+1, n+1) \in {\Bbb Z}(a,b) +  {\Bbb Z}(c,d) .
\end{eqnarray*}

In this case, we have $U = \{ 1, 2 \}$ since all of $a$, $b$, $c$ and
$d$ are positive.
By Theorem~\ref{Main} (4), we have
\begin{eqnarray*}
\mbox{$T(X; D_1, D_2)$ is a Gorenstein ring} & \Longleftrightarrow &
\mbox{$\overline{K_X+ D_1 + D_2} = 0$ in ${\rm Cl}(X)$} \\
& \Longleftrightarrow &
\mbox{$m+1 = a + c$ and $n+1 = b + d$}.
\end{eqnarray*}
\end{rm}
\end{Example}

\begin{Example}\label{ex3}
\begin{rm}
Let $a$, $b$, $c$ be pairwise coprime positive integers.
Let ${\frak p}$ be the kernel of the $k$-algebra map
$S = k[x,y,z] \rightarrow k[T]$ given by $x \mapsto T^a$,  $y \mapsto T^b$, 
$z \mapsto T^c$.

Let $\pi : X \rightarrow {\Bbb P} = {\rm Proj}(k[x,y,z])$ be the blow-up at $V_+({\frak p})$,
where $a = {\rm deg}(x)$,  $b = {\rm deg}(y)$,  $c = {\rm deg}(z)$.
Put $E = \pi^{-1}(V_+({\frak p}))$.
Let $A$ be a Weil divisor on $X$ satisfying 
$\pi^*{\mathcal O}_{\Bbb P}(1) = {\mathcal O}_X(A)$.
In this case, we have ${\rm Cl}(X) = {\Bbb Z}\overline{E} + {\Bbb Z}\overline{A} \simeq {\Bbb Z}^2$,
and $K_X = E - (a+b+c)A$.

Then, we have
\begin{eqnarray*}
{\rm Cox}(X) = R(X; -E, A) & = & R_s'({\frak p}) \mathrel{\mathop{:=}} 
S[t^{-1}, {\frak p}t, {\frak p}^{(2)}t^2, {\frak p}^{(3)}t^3, \ldots ] \subset S[t^{\pm 1}] ,\\
T(X; -E, A) & = & R_s({\frak p}) \mathrel{\mathop{:=}} 
S[{\frak p}t, {\frak p}^{(2)}t^2, {\frak p}^{(3)}t^3, \ldots ]  \subset S[t] .
\end{eqnarray*}

By Theorem~\ref{EKWHK} (4), we have
\[
\omega_{R_s'({\frak p})} = M_{K_X}= R_s'({\frak p}) (\overline{K_X}) =
R_s'({\frak p})(-1, -a-b-c) .
\]

In this case, $U = \{ 1 \}$.
By Theorem~\ref{Main} (4),
\begin{eqnarray*}
\omega_{R_s({\frak p})}  & = & M_{K_X} \cap t_1t_2 k(X)[t_1, t_2^{\pm 1}]  \\
& = & \omega_{R_s'({\frak p})} \cap t_1t_2 k(X)[t_1, t_2^{\pm 1}] \\
& = & R_s'({\frak p})(-1, -a-b-c) \cap t_1t_2 k(X)[t_1, t_2^{\pm 1}] \\
& = & R_s({\frak p}) (-1,-a-b-c) .
\end{eqnarray*}
Therefore, both of $R_s'({\frak p})$ and $R_s({\frak p})$ are quasi-Gorenstein rings, that were first proven by Simis-Trung 
(Corollary~3.4 in \cite{ST}).
Cohen-Macaulayness of such rings are deeply studied by
Goto-Nishida-Shimoda~\cite{GNS}.

Divisor class groups of ordinary and symbolic Rees rings were studied by
Shimoda~\cite{S}, Simis-Trung~\cite{ST}, etc.
\end{rm}
\end{Example}

\section{Proof of Theorem~\ref{Main}}

In this section, we shall prove Theorem~\ref{Main}.

Throughout of this section, we assume that 
$X$ is  a normal projective variety over a field $k$
such that $d = \dim X > 0$, and
$D_1$, \ldots, $D_s$ are Weil divisors on $X$ 
such that ${\mathbb N}D_1 + \cdots + {\mathbb N}D_s$ 
contains an ample Cartier divisor.


We need the following lemmta, which are well-known results.

\begin{Lemma}\label{dim}
Let $G$ be an integral domain containing a field $k$.
Let $P$ be a prime ideal of $G$.
Assume that both ${\rm tr.deg}_kG$ and ${\rm tr.deg}_kG/P$ are finite.

Then, the height of $P$ is less than or equal to
\[
{\rm tr.deg}_kG - {\rm tr.deg}_kG/P .
\]
\end{Lemma}

Using the dimension formula (e.g.\ 119p in \cite{Mat}),
it will be very easily proven.
We omit a proof.

\begin{Lemma}\label{Lemma9}
Let $r$ be a positive integer.
Let $F_1$, \ldots, $F_r$ be Weil divisors on $X$.
Let $S$ be the set of all non-zero homogeneous elements of 
$T(X; F_1, \ldots, F_r)$.
Then the following conditions are equivalent:
\begin{enumerate}
\item There exist non-negative integers $q_1$, \ldots, $q_r$ such that
$\sum_{i = 1}^rq_iF_i$ is linearly equivalent to a sum of an ample Cartier divisor and
an effective Weil divisor. 
\item There exist positive integers $q_1$, \ldots, $q_r$ such that
$\sum_{i = 1}^rq_iF_i$ is linearly equivalent to a sum of an ample Cartier divisor and
an effective Weil divisor. 
\item $S^{-1}(T(X; F_1, \ldots, F_r))
= k(X)[t_1^{\pm 1}, \ldots, t_r^{\pm 1}]$.
\item $Q(T(X; F_1, \ldots, F_r)) = 
k(X)(t_1, \ldots, t_r)$, where  $Q( \ )$ denotes the field of fractions.
\item ${\rm tr.deg}_k T(X; F_1, \ldots, F_r) = \dim X + r$.
\end{enumerate}
\end{Lemma}

It is well-known that, if $T(X; F_1, \ldots, F_r)$ is Noetherian,
then the condition (5) is equivalent to that the Krull dimension of
$T(X; F_1, \ldots, F_r)$ is $\dim X + r$.

\proof
$(2) \Rightarrow (1)$, 
and $(3) \Rightarrow (4) \Rightarrow (5)$ are trivial.

First we shall prove $(1) \Rightarrow (3)$.
Suppose
\[
\sum_{i = 1}^rq_iF_i \sim D + F ,
\]
where $q_i$'s are non-negative integers, $D$ is a very ample Cartier divisor and
$F$ is an effective divisor.
We put
\begin{eqnarray}
\label{ringC}
 C & = & \bigoplus_{m \in {\Bbb Z}}
\bigoplus_{(n_1, \ldots, n_r) \in {{\mathbb N}_0}^r}
H^0(X, {\mathcal O}_X(\sum_i n_iF_i + mD))
t_1^{n_1}\cdots t_r^{n_r}t_{r+1}^{m} \\
& \subset & k(X)[t_1, \ldots, t_r, t_{r+1}^{\pm 1}] .
\nonumber
\end{eqnarray}

We regard  $C$ as a ${\Bbb Z}^{r+1}$-graded ring with
\[
C_{(n_1, \ldots, n_r,m)} = H^0(X, {\mathcal O}_X(\sum_i n_iF_i + mD))
t_1^{n_1}\cdots t_r^{n_r}t_{r+1}^{m} .
\]
Then, we have
\[
T(X; F_1, \ldots, F_r) = \bigoplus_{(n_1, \ldots, n_r) \in {{\mathbb N}_0}^r}
C_{(n_1, \ldots, n_r,0)} ,
\]
so $T(X; F_1, \ldots, F_r)$ is a subring of $C$.
Thus, $S^{-1}C$ is a ${\Bbb Z}^{r+1}$-graded ring such that
\[
S^{-1}T(X; F_1, \ldots, F_r) = \bigoplus_{(n_1, \ldots, n_r) \in {{\mathbb N}_0}^r}
(S^{-1}C)_{(n_1, \ldots, n_r,0)} .
\]
Since $\sum_{i = 1}^rq_iF_i - D$ is linearly equivalent to 
an effective divisor $F$, there exists a non-zero element $a$ in 
\[
H^0(X, {\mathcal O}_X(\sum_i q_iF_i - D)) .
\]
For any $0 \neq b \in H^0(X, {\mathcal O}_X(D))$,
\[
(at_1^{q_1} \cdots t_r^{q_r}t_{r+1}^{-1})(bt_{r+1})
\]
is contained in $S$.
Therefore, $S^{-1}C$ contains $(bt_{r+1})^{-1}$.
Hence, $k(X)$ is 
contained in $S^{-1}C$.
Since $k(X) = (S^{-1}C)_{(0, \ldots, 0)}$, $k(X)$ is contained in $S^{-1}T(X; F_1, \ldots, F_r)$.

By the assumption of (1), 
there exists a positive integer $\ell$ such that
\[
(S^{-1}C)_{(\ell q_1, \ldots, \ell q_r, 0)} \neq 0
\]
and
\[
(S^{-1}C)_{(\ell q_1 + 1, \ell q_2, \ldots, \ell q_r, 0)} \neq 0 .
\]
Then, it is easy to see that $t_1 \in S^{-1}C$.
Therefore, $S^{-1}C$ contains $k(X)[t_1^{\pm 1}, \ldots, t_r^{\pm 1}]$.
Hence $S^{-1}T(X; F_1, \ldots, F_r)$ coincides with $k(X)[t_1^{\pm}, \ldots, t_r^{\pm}]$.

Next, we shall prove $(5) \Rightarrow (2)$.
Let $D$ be a very ample divisor.
Consider the ring
\[
R(X;F_1,\ldots,F_r,D) .
\]

First, assume that 
\[
H^0(X,{\mathcal O}_X(\sum_i u_iF_i - vD)) \neq 0
\]
for some integers $u_1$, \ldots, $u_r$, $v$
such that $v > 0$.
By the assumption (5), there exists positive integers
$u'_1$, \ldots, $u'_r$ such that 
\[
H^0(X,{\mathcal O}_X(\sum_i u'_iF_i)) \neq 0 .
\]
Therefore we may assume that there exists positive integers $u_1$, \ldots, $u_r$ and $v$
such that
\[
H^0(X,{\mathcal O}_X(\sum_i u_iF_i - vD)) \neq 0 .
\]
Here, we obtain
\[
\sum_i u_iF_i = vD + (\sum_i u_iF_i - vD) .
\]
Therefore $\sum_i u_iF_i$ is the sum of an ample divisor $vD$ 
and the divisor $\sum_i u_iF_i - vD$ which is linearly 
equivalent to an effective divisor.

Next, assume that for any integers $u_1$, \ldots, $u_r$ and $v$,
\begin{equation}\label{katei}
H^0(X,{\mathcal O}_X(\sum_i u_iF_i - vD)) = 0
\end{equation}
if $v > 0$.
We put
\[
P = 
\bigoplus_{
\superarrow{(n_1, \ldots, n_r, m) \in {\Bbb Z}^{r+1}}{m > 0}}
R(X;F_1,\ldots,F_r,D)_{(n_1, \ldots, n_r, m)} .
\]
By the assumption (5), $P$
is a prime ideal of $R(X;F_1,\ldots,F_r,D)$ of height $1$ by Lemma~\ref{dim}.
(Here, since $D$ is an ample divisor,
${\rm tr.deg}_kR(X;F_1,\ldots,F_r,D) = \dim X + r + 1$.
Remark that 
$P$ is an ideal of $R(X;F_1,\ldots,F_r,D)$ by (\ref{katei}) above.
By (5), ${\rm tr.deg}_kR(X;F_1,\ldots,F_r,D)/P = \dim X + r$.)
However $R(X;F_1,\ldots,F_r,D)$ has no homogeneous prime ideal of height $1$ 
that contains
\[
H^0(X,{\mathcal O}_X(D))t_{r+1}
\]
by Theorem~\ref{EKWHK} (2).
This is a contradiction.
\qed

Put $A = k(X)[t_1^{\pm 1}, \ldots, t_s^{\pm 1}] $ and 
$B = k(X)[t_1, \ldots, t_s]$. 
Recall that $D_1$, \ldots, $D_s$ are Weil divisors 
on a normal projective variety $X$
such that ${\mathbb N}D_1 + \cdots + {\mathbb N}D_s$ contains 
an ample Cartier divisor.
We denote $T(X; D_1, \ldots, D_s)$ and $R(X; D_1, \ldots, D_s)$
simply by $T$ and $R$, respectively.

Since
\[
T = R \cap B ,
\]
$T$ is a Krull domain.
We have proven Theorem~\ref{Main} (1).

By Theorem~\ref{EKWHK} (2), we have
\begin{eqnarray*}
R & = & \left( \bigcap_{V \in C^1(X)} R_{P_V} \right) \cap A  \\
A & = & \bigcap_{P \in {\rm NHP}^1(R)} R_{P},
\end{eqnarray*}
where ${\rm NHP}^1(R)$ is the set of non-homogeneous prime ideals of $R$ of height $1$.

It is easy to see $R_P = T_{P \cap T}$
for $P \in {\rm NHP}^1(R)$.
Therefore, we have
\[
A = \bigcap_{P \in {\rm NHP}^1(R)} T_{P \cap T} .
\]
Since $T_{P \cap T}$ is a discrete valuation ring, $P \cap T$ is a 
non-homogeneous prime ideal of $T$ of height $1$.

For $V \in C^1(X)$,  put $Q_V = P_V \cap T$.
Then, $R_{P_V} = T_{Q_V}$, since
$\sum_i{\Bbb N}D_i$ contains an ample divisor.
Therefore $Q_V$ is a homogeneous prime ideal of $T$ of height $1$.

On the other hand, 
we have $Q_i = T \cap t_i B_{(t_i)}$ and $T_{Q_i} \subset B_{(t_i)}$.
Remark that 
\[
B = A \cap ( \bigcap_{j= 1}^s B_{(t_j)} ) .
\]

Then, we have
\begin{eqnarray}
T & = & R \cap B \nonumber \\
& = & \left( \bigcap_{V \in C^1(X)} R_{P_V} \right) \cap A  \cap B \label{krull} \\
& = & \left( \bigcap_{V \in C^1(X)} T_{Q_V} \right) \cap 
\left( \bigcap_{P \in {\rm NHP}^1(R)} T_{P \cap T} \right) \cap 
\left( \bigcap_{j= 1}^s B_{(t_j)} \right) . \nonumber
\end{eqnarray}

Put 
\[
T_j = \bigoplus_{(n_1, \ldots, n_{j-1}, n_{j+1}, \ldots, n_s) \in {{\Bbb N}_0}^{s-1}}
H^0(X, {\mathcal O}_X(\sum_{i \neq j} n_iD_i))t_1^{n_1}\cdots t_{j-1}^{n_{j-1}}t_{j+1}^{n_{j+1}}\cdots t_s^{n_s} .
\]

We need the following lemma.

\begin{Lemma}\label{Le}
With notation as above, the following conditions are equivalent:
\begin{enumerate}
\item
$T_{Q_j} = B_{(t_j)}$.
\item
The height of $Q_j$ is $1$.
\item
The height of $Q_j$ is less than $2$.
\item
$j \in U$, that is, ${\rm tr.deg}_k T_j = d + s -1$.
\end{enumerate}
\end{Lemma}

\proof
By Lemma~\ref{Lemma9}, we have $Q(T) = Q(B)$.
It is easy to see that $B_{(t_j)}$ is a discrete valuation ring.
Since $Q_j$ is a non-zero prime ideal of a Krull domain $T$,
the equivalence of (1), (2) and (3) are easy.

Here, we shall prove $(1) \Rightarrow (4)$.
Remark that $T/Q_j = T_j$.
Then, we have
\[
Q(T_j) = T_{Q_j}/Q_jT_{Q_j} = B_{(t_i)}/(t_i)B_{(t_i)} = 
k(X)(t_1, \ldots, t_{j-1}, t_{j+1}, \ldots, t_s) .
\]

The implication $(4) \Rightarrow (3)$ immediately follows from
\[
{\rm ht}(Q_j) \le 
{\rm tr.deg}_k T - {\rm tr.deg}_k (T_j) = 1 .
\]
This inequality follows from Lemma~\ref{dim} and the fact
$T_j = T/Q_j$.
\qed

By (\ref{krull}), Lemma~\ref{Le} and Theorem~12.3 in \cite{Mat},
we know that 
\[
\{ Q_V \mid V \in C^1(X) \} \cup \{ Q_j \mid j \in U \}
\]
is the set of homogeneous prime ideals of $T$ of height $1$, and
\[
\{ P \cap T \mid P \in {\rm NHP}^1(R) \}
\]
is the set of non-homogeneous prime ideals of $T$ of height $1$.
Further we obtain
\[
T = \left( \bigcap_{V \in C^1(X)} T_{Q_V} \right) \cap 
\left( \bigcap_{P \in {\rm NHP}^1(R)} T_{P \cap T} \right) \cap 
\left( \bigcap_{j \in U} T_{Q_j} \right) .
\]
The proof of Theorem~\ref{Main} (2) is completed.

Let 
\[
{\rm Div}(X) = \bigoplus_{V \in C^1(X)} {\Bbb Z} \cdot V
\]
be the set of Weil divisors on $X$.
Let 
\[
{\rm HDiv}(T) = \left( \bigoplus_{V \in C^1(X)} {\Bbb Z}  \cdot  \spec(T/Q_V) \right) \oplus
\left( \bigoplus_{j \in U}  {\Bbb Z} \cdot \spec(T/Q_j) \right)
\]
be the set of homogeneous Weil divisors of $\spec(T)$.

Here, we define
\[
\phi : {\rm Div}(X) \longrightarrow {\rm HDiv}(T)
\]
by $\phi(V) = \spec(T/Q_V)$ for each $V \in C^1(X)$.
Then, it satisfies the following:
\begin{itemize}
\item
For each $a \in k(X)^\times$, we have 
\[
\phi({\rm div}_X(a)) = {\rm div}_T(a) \in 
\bigoplus_{V \in C^1(X)} {\Bbb Z}  \cdot  \spec(T/Q_V) 
\subset {\rm HDiv}(T) .
\]
\item
If $j \in U$, then 
\[
{\rm div}_T(t_j) = \spec(T/Q_j) + \phi(D_j) .
\]
\item
If $j \not\in U$, then
\[
{\rm div}_T(t_j) = \phi(D_j) .
\]
\end{itemize}
They are proven essentially in the same way as in pp631--632
in \cite{EKW}. 
Then, we have an exact sequence
\[
0 \longrightarrow \sum_{j \not\in U} {\mathbb Z}\overline{D_j}
\longrightarrow {\rm Cl}(X) \stackrel{q}{\longrightarrow}
{\rm Cl}(T)
\longrightarrow 0 
\]
such that $q(\overline{F}) = \overline{\phi(F)}$ in ${\rm Cl}(T)$.
Here, remember that ${\rm Cl}(T)$ coincides with
${\rm HDiv}(T)$ divided by homogeneous principal divisors
(e.g., Proposition~7.1 in Samuel~\cite{Sa}).

It is easy to see that the class of the Weil divisor 
$q(\overline{F})$ corresponds to the isomorphism class of
the reflexive module
\begin{eqnarray*}
& & M_F \cap \left( \bigcap_{j \in U} T_{Q_j} \right) 
= M_F \cap A \cap \left( \bigcap_{j \in U} T_{Q_j} \right) \\
& = & M_F \cap 
k(X)[t_1, \ldots, t_s, \{ t_j^{-1} \mid j \not\in U \} ]
.
\end{eqnarray*}
The proof of Theorem~\ref{Main} (3) is completed.

\begin{Remark}\label{re11}
\begin{rm}
It is easy to see 
\[
t_1^{d_1} \cdots t_s^{d_s}M_{F + \sum_id_iD_i}
= M_F 
\]
for any integers $d_1$, \ldots, $d_s$.
Therefore, we have
\begin{eqnarray*}
& & 
M_F \cap 
t_1^{d_1} \cdots t_s^{d_s}
k(X)[t_1, \ldots, t_s, \{ t_j^{-1} \mid j \not\in U \} ] \\
& = &  
t_1^{d_1} \cdots t_s^{d_s}\left(
M_{F + \sum_id_iD_i} \cap 
k(X)[t_1, \ldots, t_s, \{ t_j^{-1} \mid j \not\in U \} ]
\right) .
\end{eqnarray*}

Hence, 
\[
M_F \cap 
t_1^{d_1} \cdots t_s^{d_s}
k(X)[t_1, \ldots, t_s, \{ t_j^{-1} \mid j \not\in U \} ]
\]
is isomorphic to 
\begin{equation}\label{F +}
M_{F + \sum_id_iD_i} \cap 
k(X)[t_1, \ldots, t_s, \{ t_j^{-1} \mid j \not\in U \} ]
\end{equation}
as a $T$-module.
Remark that this is not an isomorphism as 
${\Bbb Z}^s$-graded modules.
The isomorphism class which the module (\ref{F +})
belongs to coincides with $q(\overline{{F + \sum_id_iD_i}})$.
\end{rm}
\end{Remark}

In the rest, 
we assume that $T$ is Noetherian.
We shall prove that $\omega_T$ is isomorphic to 
\[
M_{K_X} \cap 
t_1\cdots t_sk(X)[t_1, \ldots, t_s, \{ t_j^{-1} \mid j \not\in U \} ]  
\]
as a ${\Bbb Z}^s$-graded module.
(Suppose that it is true.
If we forget the grading, it is isomorphic to 
\[
M_{K_X + \sum_iD_i} \cap 
k(X)[t_1, \ldots, t_s, \{ t_j^{-1} \mid j \not\in U \} ]
\]
by Remark~\ref{re11}, 
that is corresponding to $q(\overline{K_X + \sum_iD_i})$
in ${\rm Cl}(T)$.
Therefore, we know that $\omega_T$ is $T$-free
if and only if 
\[
\overline{K_X + \sum_iD_i}
\in \sum_{j \not\in U} {\mathbb Z}\overline{D_j}
\]
in ${\rm Cl}(X)$.)

Put $X' = X \setminus {\rm Sing}(X)$.
We choose positive integers $a_1$, \ldots, $a_s$ 
and sections $f_1, \ldots, f_t \in H^0(X, \sum_ia_iD_i)$ such that
\begin{itemize}
\item
$\sum_ia_iD_i$ is an ample Cartier divisor,  
\item
$X' = \cup_k D_+ (f_k)$,  and
\item
all of the $D_i$'s are principal Cartier divisors on $D_+(f_k)$ for $k = 1, \ldots, t$.
\end{itemize}

Put $W = \{ \underline{n} \in {\Bbb Z}^s \mid \mbox{$n_i \ge 0$ if $i \in U$} \}$.
Put $D'_i = D_i|_{X'}$ for $i = 1, \ldots, s$.
Consider the morphism
\[
Y = {\rm Spec}_{X'} \left( \bigoplus_{\underline{n} \in W} 
{\mathcal O}_{X'}(\sum_in_iD'_i) t_1^{n_1} \cdots t_s^{n_s} \right)
\stackrel{\pi}{\longrightarrow} X' .
\]
Further, we have the natural map
\[
\xi : Y \longrightarrow
\spec(T) .
\]
The group ${\Bbb G}_m^s$ naturally acts on $\spec(T)$
and $Y$, and  trivially acts on $X'$.
Both $\pi$ and $\xi$ are equivariant morphisms.

\begin{Claim}\label{claim}
There exist an equivariant open subscheme $Z$ of both 
$Y$ and $\spec(T)$
such that
\begin{itemize}
\item
the codimension of $Y \setminus Z$ in $Y$ is bigger than or equal to $2$,
and
\item
the codimension of $\spec(T) \setminus Z$ in $\spec(T)$ is bigger than or equal to $2$.
\end{itemize}
\end{Claim}

\proof
For $u \in U$, 
there exist integers $c_{1u}$, \ldots, $c_{su}$ such that
\begin{itemize}
\item
$H^0(X, {\mathcal O}_X(\sum_ic_{iu}D_i)) \neq 0$,
\item
$c_{uu} = -a_u$, and
\item
$c_{iu} > 0$ if $i \neq u$.
\end{itemize}
In fact, if $u \in U$, there exist positive integers
$q_1$, \ldots, $q_{u-1}$, $q_{u+1}$, \ldots, $q_s$ such that 
\[
\sum_{i \neq u}q_iD_i
\]
is a sum of an ample divisor $D$ and a Weil divisor $F$ which is linearly equivalent to an effective divisor by Lemma~\ref{Lemma9}.
Then,
\[
H^0(X, {\mathcal O}_X(q(\sum_{i \neq u}q_iD_i) - a_uD_u) =
H^0(X, {\mathcal O}_X(q(D+F) - a_uD_u) \neq 0
\]
for $q \gg 0$.

For each $u \in U$,
we set 
\[
(b_{1u}, \ldots, b_{su}) = (c_{1u}, \ldots, c_{su}) + 
(a_1, \ldots, a_s) .
\]
Here, remark that $b_{uu} = 0$ and $b_{iu}> 0$ if $i \neq u$.

We choose 
\[
0 \neq g_u \in H^0(X, {\mathcal O}_X(\sum_ic_{iu}D_i))
\]
for each $u \in U$.

Consider the closed set of $\spec(T)$ 
defined by the ideal $J$ generated by 
\[
\{ f_kt_1^{a_1} \cdots t_s^{a_s} \mid
k = 1, \ldots, t \}
\]
and
\[
\left.
\left\{
g_u
f_k
t_1^{b_{1u}} \cdots t_s^{b_{su}} \  \right|
k = 1, \ldots, t; \ u \in U \right\} .
\]

By Theorem~\ref{Main} (2),
we know that the height of $J$ is bigger than or equal to $2$
since there is no prime ideal of $T$ of height one which contains
$J$.

We choose $d_{ki} \in k(X)^\times$ satisfying
\[
H^0(D_+(f_k), {\mathcal O}_X(D_i))
= d_{ki}H^0(D_+(f_k), {\mathcal O}_X)
\]
for each $k$ and $i$.
Then
\begin{equation}\label{Y}
Y = \bigcup_{k = 1}^t \pi^{-1}(D_+(f_k))
\ \ 
\mbox{and} \ \ 
\pi^{-1}(D_+(f_k))
= \spec(C_k) ,
\end{equation}
where
\[
C_k = H^0(D_+(f_k), {\mathcal O}_X)
[d_{k1}t_1, \ldots, d_{ks}t_s, 
\{ (d_{kj}t_j)^{-1} \mid j \not\in U \} ] .
\]

We put
\[
Z = \spec(T) \setminus V(J) .
\]
Then we have
\begin{equation}\label{Z}
Z = 
\bigcup_{k = 1}^t\left[
\spec\left( T[(f_kt_1^{a_1} \cdots t_s^{a_s})^{-1}] \right)
\cup
\left\{
\bigcup_{u \in U}
\spec\left(
T[(g_uf_kt_1^{b_{1u}} \cdots t_s^{b_{su}})^{-1}]
\right)
\right\}
\right] .
\end{equation}
Here, we have
\begin{eqnarray}\label{chart1}
T[(f_kt_1^{a_1} \cdots t_s^{a_s})^{-1}] 
& = &
H^0(D_+(f_k), {\mathcal O}_X)
[(d_{k1}t_1)^{\pm 1}, \ldots, (d_{ks}t_s)^{\pm 1}] \\
\nonumber
& = &
C_k
[ \left( \prod_{j \in U}(d_{kj}t_j)
\right)^{-1} ] .
\end{eqnarray}

On the other hand,
\begin{eqnarray*}
& & T[(g_uf_kt_1^{b_{1u}} \cdots t_s^{b_{su}})^{-1}] \\
& = & \bigoplus_{(\underline{n}) \in {\Bbb Z}^s}
T[(g_uf_kt_1^{b_{1u}} \cdots t_s^{b_{su}})^{-1}]_{(n_1, \ldots, n_s)} \\
& = & \bigoplus_{\superarrow{(\underline{n}) \in {\Bbb Z}^s}{n_u \ge 0}}
R[(g_uf_kt_1^{b_{1u}} \cdots t_s^{b_{su}})^{-1}]_{(n_1, \ldots, n_s)} \\
& = & \bigoplus_{\superarrow{(\underline{n}) \in {\Bbb Z}^s}{n_u \ge 0}}
R[
(f_kt_1^{a_1} \cdots t_s^{a_s})^{-1}, 
(g_uf_kt_1^{b_{1u}} \cdots t_s^{b_{su}})^{-1}]_{(n_1, \ldots, n_s)}
\\
& = & 
C_k[\{ (d_{kj}t_j)^{-1} \mid j \neq u \},
(g_uf_kt_1^{b_{1u}} \cdots t_s^{b_{su}})^{-1} ] .
\end{eqnarray*}
Let $\beta_{ku}$ be an element in $H^0(D_+(f_k), {\mathcal O}_X)$
such that
\[
g_u f_kt_1^{b_{1u}} \cdots t_s^{b_{su}}
= \beta_{ku} (d_{k1}t_1)^{b_{1u}} \cdots 
(d_{ks}t_s)^{b_{su}}
\]
for $k = 1, \ldots, t$ and $u \in U$.
Then, 
\begin{equation}\label{chart2}
C_k[\{ (d_{kj}t_j)^{-1} \mid j \neq u \},
(g_uf_kt_1^{b_{1u}} \cdots t_s^{b_{su}})^{-1} ]
= C_k[\left(
\beta_{ku}
\prod_{\superarrow{j \in U}{j \neq u}}
(d_{kj}t_j)
\right)^{-1}] .
\end{equation}

By (\ref{Y}), (\ref{Z}), (\ref{chart1}) and (\ref{chart2}), 
we know that $Z$ is an open subscheme of $Y$.
The ideal of $C_k$ generated by
\[
 \prod_{j \in U}(d_{kj}t_j) \ \ \mbox{and} \ \ 
\left\{ \left. \beta_{ku}
\prod_{\superarrow{j \in U}{j \neq u}}
(d_{kj}t_j) \ \right| u \in U \right\}
\]
is the unit ideal or of height two.
(If $U = \emptyset$, then $Z = Y$ by the construction.
If $U = \{ u \}$ and if $\beta_{ku}$ is a unit element,
then this ideal is the unit.
In other cases, this ideal is of height $2$.)
Therefore, 
the codimension of $Y \setminus Z$ in $Y$ is bigger than
or equal to two.
\qed

We can define the graded canonical module as in Definition~3.1 in \cite{HK}
using the theory of the equivariant twisted inverse functor~\cite{H}.

By Claim~\ref{claim} above and Remark~3.2 in \cite{HK},
we have $\omega_T = H^0(Y, \omega_Y)$.
On the other hand, we have
\begin{eqnarray*}
\omega_Y & = & \bigwedge^s \Omega_{Y/X'} \otimes \pi^*{\mathcal O}_{X'}(K_{X'})  \\
& = & \pi^*{\mathcal O}_{X'}(\sum_iD'_i)(-1, \ldots, -1)
\otimes_{{\mathcal O}_Y} \pi^*{\mathcal O}_{X'}(K_{X'}) \\
& = & \pi^*{\mathcal O}_{X'}(\sum_iD'_i + K_{X'})(-1, \ldots, -1) ,
\end{eqnarray*}
where $(-1, \ldots, -1)$ denotes the shift of degree
(Theorem~28.11 in \cite{H}).

Then, we have
\[
H^0(Y, \omega_Y) = 
H^0(X', \pi_*\pi^*{\mathcal O}_{X'}(\sum_iD'_i + K_{X'})(-1, \ldots, -1)) .
\]
By the projection formula (Lemma~26.4 in \cite{H}),
\begin{eqnarray*}
& & \pi_*\pi^*{\mathcal O}_{X'}(\sum_iD'_i + K_{X'})(-1, \ldots, -1) \\
& = & \left( {\mathcal O}_{X'}(\sum_iD'_i + K_{X'}) \otimes 
\pi_*{\mathcal O}_{Y} \right) (-1, \ldots, -1) \\
& = &  \left( {\mathcal O}_{X'}(\sum_iD'_i + K_{X'}) \otimes 
\left[
\bigoplus_{\underline{n} \in W} {\mathcal O}_{X'}(\sum_in_iD'_i)
\right] \right) (-1, \ldots, -1) \\
& = & 
\left(
\bigoplus_{\underline{n} \in W} {\mathcal O}_{X'}(\sum_i(n_i+1)D'_i
+ K_{X'})
\right) (-1, \ldots, -1) \\
& = & 
\bigoplus_{\underline{n} \in W+(1,\ldots,1)} 
{\mathcal O}_{X'}(\sum_in_iD'_i+ K_{X'}) .
\end{eqnarray*}

Therefore, we have
\begin{eqnarray*}
H^0(Y, \omega_Y) & = &
H^0(X', \bigoplus_{\underline{n} \in W+(1,\ldots,1)} 
{\mathcal O}_{X'}(\sum_in_iD'_i+ K_{X'}) ) \\
& = & 
\bigoplus_{\underline{n} \in W+(1,\ldots,1)} 
H^0(X', {\mathcal O}_{X'}(\sum_in_iD'_i+ K_{X'}) ) \\
& = & 
\bigoplus_{\underline{n} \in W+(1,\ldots,1)} 
H^0(X, {\mathcal O}_{X}(\sum_in_iD_i+ K_{X}) ) \\
& = & 
M_{K_X} \cap t_1 \cdots t_sk(X)[t_1, \ldots, t_s, \{ t_j^{-1} \mid j \not\in U \} ] .
\end{eqnarray*}

We have completed the proof of Theorem~\ref{Main}.

\noindent
Department of Mathematics \\
School of Science and Technology \\
Meiji University \\
Higashimita 1-1-1, Tama-ku \\
Kawasaki 214-8571, Japan \\

\vspace{1mm}

\noindent
{\tt kurano@isc.meiji.ac.jp} \\
{\tt http://www.math.meiji.ac.jp/\~{}kurano}


\begin{thebibliography}{99}
\bibitem{EKW}
{\sc E. Javier Elizondo, K. Kurano and K.-i. Watanabe},
{\em The total coordinate ring of a normal projective variety},
J. Algebra {\bf 276} (2004), 625--637.

\bibitem{GNS}
{\sc S. Goto, K. Nishida and Y. Shimoda},
{\em The Gorensteinness of symbolic Rees algebras for space curves},
J. Math. Soc. Japan {\bf 43} (1991), 465–-481.

\bibitem{H}
{\sc M. Hashimoto},
{\em Equivariant Twisted Inverses}, in 
{\em Foundations of Grothendieck Duality for Diagrams of Schemes}
(J. Lipman, M. Hashimoto, eds.),
Lecture Notes in Math. {\bf 1960}, Springer (2009), pp.~261--478.

\bibitem{HK} {\sc M. Hashimoto and  K. Kurano},
{\it The canonical module of a Cox ring},
Kyoto J. Math. {\bf 51} (2011), 855--874.

\bibitem{HuKeel} {\sc Y. Hu and  S. Keel},
{\it Mori dream spaces and GIT},
Michigan Math J. {\bf 48} (2000), 331--348.

\bibitem{Mat} {\sc H. Matsumura},
{\it Commutative ring theory},
Cambridge University Press, 1990.

\bibitem{M} {\sc S. Mori},
{\it Graded factorial domains},
Japan J. Math. {\bf 2} (1977), 223--237.

\bibitem{Sa} {\sc P. Samuel},
{\it Lectures on unique factorization domains}, 
Tata Inst. Fund. Res., Bombay, 1964.

\bibitem{S} {\sc Y. Shimoda},
{\it The class group of the Rees algebras over polynomial rings},
Tokyo J. Math. {\it 2} (1979), 129--132.

\bibitem{ST} {\sc A. Simis and N. V. Trung},
{\it The divisor class group of ordinary and symbolic blow-ups},
Math. Z. {\bf 198} (1988), 479--491.

\bibitem{W}
{\sc K.-i. Watanabe},
{\em Some remarks concerning Demazure's construction of normal graded rings},
Nagoya Math. J. {\bf 83} (1981),  203--211.
\end{thebibliography}
\end{document}